\documentclass[reqno]{amsart}
\usepackage{amsmath,amssymb,amsthm}
\usepackage{bm}
\usepackage{extarrows}
\usepackage{mathrsfs}
\usepackage{stmaryrd}
\usepackage[colorlinks=true]{hyperref}
\usepackage[all]{xy}

\newtheorem{theorem}{Theorem}[section]

\newtheorem{proposition}[theorem]{Proposition}

\theoremstyle{definition}
\newtheorem{definition}[theorem]{Definition}
\newtheorem{example}[theorem]{Example}

\theoremstyle{remark}

\numberwithin{equation}{section}

\newcommand\floor[1]{{\left[#1\right]}}
\newcommand\surj{\twoheadrightarrow}
\newcommand\set[1]{\left\{#1\right\}}
\newcommand\bfu{{\mathbf{u}}}
\newcommand\bfv{{\mathbf{v}}}
\newcommand\bfN{{\mathbf{N}}}
\newcommand\BC{{\mathbb{C}}}
\newcommand\BF{{\mathbb{F}}}
\newcommand\BN{{\mathbb{N}}}
\newcommand\BQ{{\mathbb{Q}}}
\newcommand\BR{{\mathbb{R}}}
\newcommand\BZ{{\mathbb{Z}}}
\newcommand\CO{{\mathcal{O}}}
\newcommand\fl{{\mathfrak{l}}}
\newcommand\diag{{\mathrm{diag}}}
\newcommand\Gal{{\mathrm{Gal}}}
\newcommand\GL{{\mathrm{GL}}}
\newcommand\lcm{\mathrm{lcm}}
\newcommand\ord{{\mathrm{ord}}}
\newcommand\Tr{{\mathrm{Tr}}}
\newcommand\Vol{{\mathrm{Vol}}}

\title[The virtual periods of LRS]{The virtual periods of linear recurrence sequences in cyclotomic fields}
\author{Shenxing Zhang}
\address{School of Mathematics, Hefei University of Technology, Hefei, Anhui 230009, China}
\email{zsxqq@mail.ustc.edu.cn}

\date{\today}

\begin{document}
\maketitle
\begin{abstract}
A linear recurrence sequence in a cyclotomic field produces a sequence of the generating fields of each term. We show that the later sequence is periodic after removing the first finite terms, and give a bound of its period. This can be applied to exponential sums.
\end{abstract}
\tableofcontents

\section{The virtual periods of linear recurrence sequences}

A {\em linear recurrence sequence} with dimension $n$ is a sequence $\set{a_k}_{k\ge 0}$ such that for all $k\ge n$, 
  \[a_k=\sum_{i=1}^n c_i a_{k-i}\]
for some constant $c_1,\dots,c_n$ where $c_n\neq 0$.

\begin{definition}
We say a sequence $\set{a_k}_{k\ge 0}$ is \emph{virtually periodic} if there exists integers $N,r\ge 1$ such that $a_n=a_{n+r}$ for any $n\ge N$. The minimal $r$ is called the \emph{virtual period} and the minimal $N$ is called the \emph{pre-period length}.
\end{definition}

One can easily obtain the following result from \cite[Theorem~3]{WanYin2020}.
\begin{theorem}
Let $K$ be a field of characteristic $0$ and $L$ a finite extension of $K$.
Let $\set{a_k}_{k\ge 0}$ be a linear recurrence sequence in $L$. Then the sequence $\set{K(a_k)}_{k\ge 0}$ of fields is virtually periodic.
\end{theorem}

For $L=\BQ(\mu_m)$ a cyclotomic field, we will reprove the Skolem-Mahler-Lech Theorem and the theorem above to estimate the virtual period, see \cite{Hansel1985}, \cite[\S 2]{Lech1953}.

Certainly, we may assume that $2\mid m$.
For $n>1$, denote by
\begin{equation}\label{eq:ep}
	e_p=\begin{cases}
		0, &p\nmid m, p>n+1;\\
		1+\floor{\log_p\frac{n}{p-1}}, &p\nmid m, 2<p\le n+1;\\
		2+\floor{\log_2 n}, &p=2,\ord_2(m)=1;\\
		\ord_p(m)+\floor{\log_p n}, &2p\mid m
	\end{cases}
\end{equation}
for each prime $p$, where $\floor{x}$ denote the greatest integer less than or equal $x$.
Denote by
\begin{equation}\label{eq:Rmn}
	R_{m,n}=\prod_p p^{e_p}
\end{equation}
and $R_{m,1}=m$. For odd $m$, we denote $R_{m,n}=R_{2m,n}$.

\begin{theorem} \label{thm:virtual_period}
Let $\set{a_k}_{k\ge 0}$ be a linear recurrence sequence in $\BQ(\mu_m)$ with dimension $n$.

(1) There exists a positive integer $s\mid R_{m,n}$ such that the set $\set{k: a_k=0}$ is a union of some $i+s\BN$ and a finite set.

(2) The sequence $\set{\BQ(a_k)}_{k\ge 0}$ is virtually periodic of virtual period $r\mid R_{m,n}$.
Moreover, $\BQ(a_k)\subseteq \BQ(a_{k'})$ if $k\equiv k'\bmod r$ and $k'\ge N$ for pre-period length $N$.
\end{theorem}

\begin{proof}
(1) Let $\lambda$ be a positive integer which is a common multiplier of the denominators of $a_0,\dots,a_{n-1}, c_1,\dots,c_n$.
Then $a_k':=a_k \lambda^{k+1}$ satisfies
  \[a_k'=\sum_{i=1}^n c_i \lambda^i a_{k-i}'\]
and $a_0',\dots,a_{n-1}',c_1\lambda,\dots,c_n\lambda^n \in \BZ[\mu_m]$.
Thus we may assume that $c_1,\dots,c_n$ and all $a_k$ lie in $\BZ[\mu_m]$. 

Let $M$ be the $n\times n$ matrix with $M_{i,1}=c_i,M_{i,i+1}=1$ for all $i$, and other entries are all zero.
Denote
  \[\bfu=(a_{n-1},a_{n-2},\dots,a_0),\quad \bfv=(1,0,\dots,0)^T.\]
Then $a_k=\bfu M^{k+1-n}\bfv.$

Let $\ell>2$ be a prime splits completely in $K=\BQ(\mu_m)$, such that $\ell\nmid c_n=(-1)^{n-1}\det M$.
Let $\fl$ be a prime of $K$ above $\ell$ and denote by $\CO_{\fl}$ the completion of $\BZ[\mu_m]$ at $\fl$.
Then $\ell$ is a uniformizer of $\CO_\fl$ and the residue field is $\kappa_\fl\cong\BF_\ell$.
Denote by $s(\ell)$ the order of the image of $M$ under $\CO_{\fl}\surj \kappa_\fl$.
Then $M^{s(\ell)}=I+\ell M'$ for some matrix $M'$ over $\CO_\fl$.
For $i\ge n-1$, the function
  \[a_{i+s(\ell)x}:=\bfu M^{i+1-n}(I+\ell M')^x\bfv=\sum_{k\ge 0}{x\choose k}\bfu M^{i+1-n}M'^k \bfv\cdot \ell^k\]
on $x\in\CO_\fl$ converges since $\ord_\ell\bigl(\ell^k/k!\bigr)> k\frac{\ell-2}{\ell-1}$ tends to infinity.
If there are infinitely many integers $x$ such that $a_{i+s(\ell)x}=0$, then the set of these indices has an accumulation point in $\CO_\fl$ in $\ell$-adic topology.
Hence the function $a_{i+sx}$ must be zero identically by Weierstrass preparation theorem.
From this we know that the set $\set{k: a_k=0}$ has the predicated form.

If $s$ is a positive integer such that $\set{k: a_k=0}$ is a union of some $i+s\BN$ and a finite set, then so is its multipliers. We take the minimal $s$. Then $s\mid s(\ell)$ is the order of an element in $\GL_n(\kappa_\fl)=\GL_n(\BF_\ell)$. 
We will use the following proposition to estimate $s$.

\begin{proposition}[{\cite[\S 1, Corollary~1]{Darafsheh2005}}]\label{pro:max_order}
Each maximal order of elements in $\GL_n(\BF_\ell)$ has form
	\[\ell^t\times \lcm(\ell^{d_1}-1,\dots,\ell^{d_s}-1),\]
where $\sum_{i=1}^s k_i d_i=n$ has integer solutions and $t$ is the smallest non-negative integer such that $\ell^t\ge \max\set{k_1,\dots,k_s}$.
In particular, the $p$-order of each maximal order is at most $\max_{d\le n}\ord_p(\ell^d-1)$.
\end{proposition}

We may assume that $n\ge 2$.
For any rational prime $p\nmid m$, there exists a prime $\ell\nmid c_n$ which splits completely in $K$ and $\ell$ is a primitive root modulo $p^2$.
Thus
	\[\ord_p(s)\le \max_{d\le n} \ord_p(\ell^d-1)=\begin{cases}
		0,&p>n+1;\\
		1+\floor{\log_p\frac{n}{p-1}},\ &2<p\le n+1.
	\end{cases}\]

There exists a prime $\ell\nmid c_n$ which splits completely in $K$ and $\ell\not\equiv 1\bmod pm$ for any $p\mid m$.
Then for each $p\mid m$, we have
  \[\ord_p(s)\le\max_{d\le n}\ord_p(\ell^d-1)=
		\begin{cases}
		2+\floor{\log_2 n},&p=2,\ord_2(m)=1;\\
		\ord_p(m)+\floor{\log_p n},&\text{otherwise}.
	\end{cases}\]

(2) For $\sigma\in\Gal\bigl(\BQ(\mu_m)/\BQ)\bigr)$,	 $\set{\sigma(a_k)-a_k}$ is a linear recurrence sequence satisfying
	\[\sigma(a_k)-a_k=(\sigma\bfu,\bfu)\left(\begin{smallmatrix}\sigma M&\\&M	\end{smallmatrix}\right)^{k+1-n}
		\left(\begin{smallmatrix}\bfv\\-\bfv\end{smallmatrix}\right).\]
Similar to (1), let $\ell>2$ be a prime splits completely in $K$ such that $\ell\nmid c_n c_n^\sigma$.
Then 
	\[M^{s(\ell)}\equiv (M^\sigma)^{s'(\ell)}\equiv I\bmod\fl,\]
where $s(\ell),s'(\ell)$ are orders of two elements in $\GL_n(\BF_\ell)$. Thus the set $\set{k: \sigma(a_k)=a_k}$ is a union of a finite set and some $i+r_\sigma\BN$, where $r_\sigma\mid\lcm\bigl(s(\ell),s'(\ell)\bigr)$.
By Proposition~\ref{pro:max_order} and the estimation in (1), we have $\ord_p(r_\sigma)\le e_p$ for each prime $p$.

Denote by $r$ the least common multiplier of these $r_\sigma$. Then there exists $N$ such that $\sigma(a_k)=a_k$ if and only if $\sigma(a_{k+r})=a_{k+r}$ for any $k\ge N$. 
Denote by $H_k$ the set of $\sigma\in\Gal(K/\BQ)$ fixing $a_k$.
As shown in \cite[Theorem~3]{WanYin2020}, $H_k=H_{k+r}$ for any $k>N$. Hence $\BQ(a_k)=\BQ(a_{k+r})$. Certainly, $r\mid R_{m,n}$.

For any integer $k\ge 0$, denote by $k'$ the minimal one such that $k_0\equiv k\bmod r$ and $k'\ge N$.
Then $\sigma\in H_{k'}$ fixes $a_{k'+ir}$ for any $i\ge0$ and the sequence $\set{\sigma(a_{k'+ir})-a_{k'+ir}}_{i\ge0}$ is identically zero.
This implies that $\set{\sigma(a_{k+ir})-a_{k+ir}}_{i\ge 0}$ is identically zero since it is a linear recurrence sequence.
Hence $a_k$ is fixed by $H_{k'}$ and $a_k\in\BQ(a_{k'})$.
\end{proof}

\section{The virtual periods of exponential sums}

Let $f\in\BF_q[x_1^{\pm1},\dots,x_m^{\pm1}]$ be a Laurent polynomial and $\chi:(\BF_q^\times)^m\to \BC^\times$ a character of order $c$.
Define the \emph{(toric) exponential sums}
	\[S_k(f,\chi)=\sum_{x \in(\BF_{q^k}^\times)^m} \psi\Bigl(\Tr_{\BF_{q^k}/\BF_p}\bigl(f(x)\bigr)\Bigr)\chi\bigl(\bfN_{\BF_{q^k}/\BF_q}(x)\bigr)\in\BZ[\mu_{pc}].\]
Then the $L$-function
  \[L(T,f,\chi)=\exp\left(\sum_{k=1}^\infty \frac{T^k}{k}S_k(f,\chi)\right)\]
is a rational function over $\BQ(\zeta_{pc})$ by the Dwork-Bombieri-Grothendick rationality theorem (\cite{Bombieri1966}).
Write
	\[L(T,f,\chi)=\frac{\prod_j(1-\beta_j T)}{\prod_i (1-\alpha_i T)}.\]
Then
	\[S_k(f,\chi)=\sum_i \alpha_i^k-\sum_j \beta_j^k\]
and $\set{S_k(f,\chi)}_{k\ge 1}$ is a linear recurrence sequence in $\BQ(\mu_{pc})$. Hence we have:

\begin{theorem} \label{co:virtual_period}
The sequence $\set{\BQ\bigl(S_k(f,\chi)\bigr)}_{k\ge 1}$ is virtually periodic with the period dividing $R_{pc,n}$, where $n$ is the number of zeroes and poles of the $L$-function $L(T,f,\chi)$ and $c$ is the order of $\chi$.
In particular, every prime factors of the virtual period are less than $n+1$ or divides $pc$.
\end{theorem}

We will omit $\chi$ if it's trivial.
\begin{example}
Assume that $d$ is a divisor of $q-1$.
Let $f(x)=x^d+a\in\BF_q[x]$ be a polynomial and $\chi$ the trivial character.
Then by the Hasse-Davenport relation, we have
	\[S_k(f)=-\sum_{i=1}^{d-1} \beta_i^k,\quad L(T,f)=\prod_{i=1}^{d-1}(1-\beta_i T),\]
where
  \[\beta_i=-\psi\bigl(\Tr_{\BF_q/\BF_p}(a)\bigr)\tau(\omega^{\frac{(q-1)i}d})\in\BQ(\mu_{pd}),\quad \tau(\eta)=\sum_{x\in\BF_q^\times}\eta(x)\psi\bigl(\Tr_{\BF_q/\BF_p}(x)\bigr).\]
Then
	\[a_k=\bfu M^k\bfv,\quad \bfu=-\bfv^T=(1,\dots,1),\  M=\diag\set{\beta_1,\dots,\beta_{d-1}}.\]
Similarly to the proof of Theorem~\ref{thm:virtual_period}, we take a prime $\ell>2$ which splits completely in $\BQ(\mu_{pd})$ and $\ell\nmid \det M$. Then $M^{\ell-1}\equiv I\bmod\ell$.
Hence the virtual period of $\set{\BQ\bigl(S_k(f)\bigr)}_k$ divides $R_{pd,1}=pd$ or $2pd$.

If $d|(p-1)$ or $d\mid\frac{q-1}{p-1}$, we have $\deg S_k(x^d)=\frac{p-1}{\bigl(p-1,(q^k-1)/d\bigr)}$ by \cite[Theorem~4.8]{Wan2021} and the degree sequence $\set{\BQ\bigl(S_k(f)\bigr)}_k$ is periodic with the period 
	\[\begin{cases}
		d,&\text{if}\ d\mid(p-1), \ \Tr_{\BF_q/\BF_p}(a)=0;\\
		pd,&\text{if}\ d\mid(p-1), \ \Tr_{\BF_q/\BF_p}(a)\neq 0;\\
		1,&\text{if}\ d\mid\frac{q-1}{p-1}, \ \Tr_{\BF_q/\BF_p}(a)=0;\\
		p,&\text{if}\ d\mid\frac{q-1}{p-1}, \ \Tr_{\BF_q/\BF_p}(a)\neq 0.
	\end{cases}\]
\end{example}

\begin{example}
The exponential sums $S_k(f)$ of $f=x_1+\cdots+x_{n-1}+\frac{a}{x_1\cdots x_{n-1}}$ are called \emph{Kloosterman sums}.
The sequence $\set{\BQ\bigl(S_k(f)\bigr)}_{k}$ is virtually periodic with the period dividing $R_{p,n}$, since $L(T,f)^{(-1)^n}$ is a polynomial of degree $n$ by Deligne in \cite{Deligne1977}.

It's known that if $p$ is large with respect to $\log_p q, n$ and $c$, then the generating fields of twisted Kloosterman sums are known. See \cite{Fisher1992} and \cite{Zhang2021a}.
\end{example}

\begin{example}
It's easy to see that $R_{m,n}\mid R_{m,n+1}$.
Bombieri in \cite[Theorem~1]{Bombieri1978} showed that the number of zeroes and poles of $S_k(f)$ is at most $4d+5$, where $d$ is the degree of $f$ (be careful the different definitions of exponential sums).
Hence the virtual period of $\set{\BQ\bigl(S_k(f)\bigr)}_{k\ge 1}$ divides $R_{p,4d+5}$.
\end{example}

\begin{example}
If $f$ is so-called \emph{non-degenerate}, then $L(T,f,\chi)^{(-1)^{n-1}}$ is a polynomial of degree $n!\Vol(\Delta)$, where $\Delta$ is the convex polyhedron in $\BR^n$ associated to $f$. 
Hence the virtual period of $\set{\BQ\bigl(S_k(f,\chi)\bigr)}_{k\ge 1}$ divides $R_{p,n!\Vol(\Delta)}$.
See \cite[Corollary 2.12]{AdolphsonSperber1993}, \cite[Theorem~1.3]{LiuWei2007} and \cite[Theorem~1]{Liu2007}.

In particular, if $f$ is a polynomial in one variable with $p\nmid d=\deg f$, then the virtual period of $\set{\BQ\bigl(S_k(f,\chi)\bigr)}_{k\ge 1}$ divides $R_{p,d}$.
\end{example}

\textbf{Acknowledgments.}
The author would like to thank Daqin Wan and Yi Ouyang for many helpful suggestions.
This work is partially supported by NSFC (Grant No. 12001510), the Fundamental Research Funds for the Central Universities (No. WK0010000061) and Anhui Initiative in Quantum Information Technologies (Grant No. AHY150200).

\end{document}